\documentclass[%
reprint,
amsmath,amssymb,
aps,
]{revtex4-2}

\usepackage[normalem]{ulem}
\usepackage{graphicx}% Include figure files
\usepackage{dcolumn}% Align table columns on decimal point
\usepackage{bm}% bold math
\usepackage{amssymb,amsmath,amsfonts}
\usepackage{color}

\usepackage[dvipsnames]{xcolor}
\usepackage{lineno}
\newcommand{\mbX}{\mathbf{X}}
\newcommand{\mbY}{\mathbf{Y}}
\newcommand{\mbx}{\mathbf{x}}
\newcommand{\mby}{\mathbf{y}}
\newcommand{\mbf}{\mathbf{f}}
\newcommand{\mbF}{\mathbf{F}}
\newcommand{\mbg}{\mathbf{g}}
\newcommand{\E}{\mathbb{E}}
\newcommand{\R}{\mathbb{R}}
\newcommand{\LL}{\mathcal{L}}
\newcommand{\LLd}{\mathcal{L}^\dagger}
\def\given{\:|\:}

\begin{document}
	
	\preprint{APS/123-QED}
	
	\title{Isostables for stochastic oscillators}% Force line breaks with \\
	
	\author{Alberto P\'erez-Cervera}
	% \altaffiliation[Also at ]{}%Lines break automatically or can be forced with \\
	%\author{}%
	%\email{Second.Author@institution.edu}
	%\affiliation{%
		% Authors' institution and/or address\\
		% This line break forced with \textbackslash\textbackslash
		%}%
	
	%\collaboration{MUSO Collaboration}%\noaffiliation

	\affiliation{
		%Center for Cognition and Decision Making, Institute for Cognitive Neuroscience, 
		National Research University Higher School of Economics, Moscow, Russia.\\
		\& Instituto de Matemática Interdisciplinar, Universidad Complutense de Madrid, 28040 Madrid, Spain
	}%
	\author{Benjamin Lindner}
	\affiliation{%
			Bernstein Center for Computational Neuroscience Berlin, Berlin, Germany\\ 
		\& 	Department of Physics, Humboldt University, Berlin, Germany.
	}%
	\author{Peter J. Thomas}
	\affiliation{%
	Department of Mathematics, Applied Mathematics, and Statistics, Case Western Reserve University, Cleveland, Ohio 44106, USA
	}%

	\date{\today}% It is always \today, today,
	%  but any date may be explicitly specified
	
	\begin{abstract}
		%While the analysis of deterministic oscillators through phase-amplitude reduction is well understood, an analogous framework for stochastic oscillators remains elusive. 
		%We extend the notion of the isostable coordinate for stochastic oscillators via the eigenfunction with the smallest magnitude nonzero real eigenvalue of the adjoint Kolmogorov operator. 
		%Together with the stochastic asymptotic phase (defined via the eigenfunctions with complex eigenvalues having least negative real part [Thomas and Lindner, Phys.~Rev.~Lett. (2014)]) the stochastic isostable function completes a phase-amplitude description of stochastic oscillators. 
		%This framework provides a physically well-founded definition of a \emph{stochastic limit cycle}, both for deterministic systems perturbed by noise, and for systems with noise-induced oscillations, such as stochastic heteroclinic oscillators and noisy spiral foci. (894 char)
		
		%\textcolor{blue}{In this letter we complete the phase-amplitude description of stochastic oscillators by extending the isostable coordinate to the stochastic domain. While the eigenfunctions of the adjoint Kolmogorov operator with the slowest decaying complex eigenvalues defined a stochastic asymptotic phase, the stochastic isostables correspond to the eigenfunction having the smallest magnitude nonzero real eigenvalue. Our framework provides a physically well-founded definition of a stochastic limit cycle, even for noise-induced oscillators, such as stochastic heteroclinic oscillators or noisy spiral foci. (597 char)}
		
		Thomas and Lindner (2014, Phys.~Rev.~Lett.)~defined an asymptotic phase for stochastic oscillators as the angle in the complex plane made by the eigenfunction, having complex eigenvalue with least negative real part, of the backward Kolmogorov (or stochastic Koopman) operator.  
		We complete the phase-amplitude description of noisy oscillators by defining the stochastic isostable coordinate as the eigenfunction with the least negative nontrivial real eigenvalue. 
		%(sub1:) This framework provides a physically well-founded definition of ``stochastic limit cycle" that includes noise-induced oscillations. 
		%\rot{Our results suggest a definition of ``stochastic limit cycle" that encompasses noise-induced oscillations.}	
		Our results suggest a framework for stochastic limit cycle dynamics that encompasses noise-induced oscillations.
		
		%Phase-amplitude reduction for a stable, deterministic limit cycle oscillator seeks a change of coordinates giving a phase variable that evolves at a constant rate, and one or more transverse "isostable" coordinate variables that decay according to linear stable ordinary differential equations. Several generalizations of phase reduction to the case of stochastic oscillatory systems have been proposed, including introduction of "isochrons" with a constant mean return time property, and extraction of an "asymptotic phase" from the complex angle function extracted from the eigenfunction associated with the slowest decaying nontrivial complex eigenvalue of the generator of the underlying Markov process (also known as the Koopman operator or the adjoint Kolmogorov operator). Here we describe a generalization of the isostable coordinate to an oscillatory stochastic system. Our stochastic isostable is derived from the eigenfunction associated with the slowest decaying nontrivial real eigenvalue.
		%\begin{description}
			%\item[Usage]
			%Secondary publications and information retrieval purposes.
			%\item[Structure]
			%You may use the \texttt{description} environment to structure your abstract;
			%use the optional argument of the \verb+\item+ command to give the category of each item. 
			%\end{description}
	\end{abstract}
	
	%\keywords{Suggested keywords}%Use showkeys class option if keyword
	%display desired
	\maketitle
	
	%\tableofcontents
	
	%\section{\label{sec:level1}First-level heading:\protect\\ The line
		%break was forced \lowercase{via} \textbackslash\textbackslash}
	
	\textit{Introduction -} Nonlinear stochastic oscillations occur throughout natural and engineered systems.  
	Examples include the membrane potential of nerve cells \citep{BryantMarcosSegundo1973JNP,WalterEtAlKhodakhah2006NatNeuro}, 
	the deflection of sensory organelles \cite{MarBoz03}, the concentration of intracellular calcium \cite{SkuKet08}, 
	the populations of predators and prey \cite{mckane2005predator}, 
	chemical reactions \cite{feistel1978deterministic}, climate systems \cite{GanRah02} 
	and the intensity of lasers \citep{Mcnamara1988observation}.
	For a system described by deterministic 
	ordinary differential equations,
	\begin{equation}
		\label{eq:ODE}
		\frac{d\mbx}{dt}=\mbf(\mbx),\quad \mbx\in\R^n
	\end{equation}
	robust oscillations arise from orbitally stable limit cycle (LC) solutions $\mbx=\gamma(t)=\gamma(t+T)$ with finite period $T$.
	The analysis of limit cycle systems through phase-amplitude reduction provides an essentially complete understanding of oscillatory dynamics \cite{guillamon2009computational, WilsonErmentrout2019PRL-phase}.    
	By transforming to
	a phase variable $\theta$, together with $n-1$ amplitude variables $\sigma$ that decay at  rates given by the non-trivial Floquet exponents $\lambda_i$, ($i=1,\dots,n-1$), the system eq.~\eqref{eq:ODE} becomes
	\begin{equation}
		\label{eq:phase-amplitude-deterministic}
		\dot\theta=2\pi/T,\qquad\dot\sigma=\Lambda\sigma 
	\end{equation}
	with $\Lambda = \text{diag}(\lambda_1, \dots, \lambda_{n-1})$ \citep{ Wilson2020,perez2020global}.  
	This transformation facilitates the understanding of synchronization, entrainment and control in a broad range of scenarios \cite{pikovsky2003synchronization,castejon2013phase, shirasaka2017phase,wilson2018greater,MongaWilsonMatchenMoehlis2019BICY-phase}. Indeed, whereas the level sets of $\theta$ form the well known ``isochrons" \cite{guckenheimer1975}, the LC coincides with the set of points $\{\mbx\given\sigma(\mbx)=0\}$.
	
	In contrast, in stochastic systems, oscillatory behavior can arise not only from an underlying limit cycle, but in a noise-perturbed spiral--fixed-point system (quasicycle) (see the early example \cite{UhlOrn30}), in systems with an underlying stable heteroclinic orbit driven by random fluctuations \cite{Giner-Baldo2017JStatPhys}, or in noisy excitable systems \cite{lindner2004effects}. 
	Incorporating stochastic dynamics fundamentally changes the conceptual foundations underlying the phase-amplitude construction.

	Consider the Langevin equation 
	\begin{equation}
		\label{eq:SDE}
		\frac{d\mbX}{dt}=\mbf(\mbX) + \mbg(\mbX)\xi(t)
	\end{equation} 
	where $\mbf$ is an $n$-dimensional vector,  $\mbg$ is an $n\times k$ matrix, $\xi$ is $k$-dimensional white noise with uncorrelated components \footnote{For an example in which $k\not=n$, see \cite{pu2021resolving,pu2020fast}}, $\langle\xi_i(t)\xi_j(t')\rangle=\delta(t-t')\delta_{i,j}$, and we interpret the stochastic differential equation in the sense of Itô \footnote{We choose the Itô interpretation for its  mathematical convenience. 
		For every Stratonovich-interpreted SDE there is an equivalent Itô-interpreted SDE \cite{gardiner1985handbook}.  
		Thus, choosing between the Itô or the Stratonovich interpretation will not change our framework, which is based on the (uniquely defined) backward Kolmogorov operator.}.
	For this system, orbits are no longer periodic. 
	The transit times between isochrons or other Poincar\'e sections are random variables \cite{schwabedal2013phase, cao2020partial}. 
	%The phase and amplitude variables eq.~\eqref{eq:phase-amplitude-deterministic} evolve stochastically.
	The classical notion of a ``limit cycle", as a closed, isolated periodic orbit, is no longer well defined, and a shift in perspective is required.
	
	By considering the evolution of ensembles of trajectories, two of us \cite{thomas2014asymptotic} established a generalization of the asymptotic phase to stochastic oscillators, in terms of the slowest decaying complex eigenfunction of the generator $\LLd$ of the Markov process eq.~\eqref{eq:SDE} \footnote{For an alternative approach to phase reduction for stochastic oscillators see \cite{schwabedal2013phase,cao2020partial}.}.  
	In this Letter we generalize the amplitude coordinate for a planar stochastic oscillator, which we define to be the slowest decaying \emph{real} eigenfunction of $\LLd$.
	%By putting the notion of phase and amplitude on a common basis, we complete the foundation for the study of stochastic oscillatory systems, and provide a physically grounded definition for a ``stochastic limit cycle".
	By putting the notion of phase and amplitude on a common basis, we suggest a framework for the study of stochastic oscillatory systems.
	
	\textit{Mathematical Preliminaries -} As in \cite{thomas2014asymptotic,thomas2019phase}, we describe an ensemble of trajectories through the density \small
	\begin{equation*}
		\rho(\mby,t\given\mbx,s)=\frac{1}{|d\mby|}\Pr\left\{\mbX(t)\in[\mby,\mby+d\mby)\given \mbX(s)=\mbx \right\}
	\end{equation*} \normalsize
	for $s<t$.
	The density satisfies Kolmogorov's  equations \small
	\begin{align*}
		\frac{\partial}{\partial t}\rho(\mby,t\given\mbx,s)&=\LL_\mby[\rho]=-\nabla_\mby\cdot\left( \mbf(\mby) \rho \right)+\frac{\partial^2}{\partial y_i y_j}\left(D_{ij}(\mby)\rho\right)\\
		-\frac{\partial}{\partial s}\rho(\mby,t\given\mbx,s)&=\LLd_\mbx[\rho]=\mbf(\mbx)\cdot\nabla_\mbx\left( \rho \right)+D_{ij}(\mbx)\frac{\partial^2}{\partial x_i x_j}\left(\rho\right)
	\end{align*}
	\normalsize
	where $D=\frac12 gg^\intercal$.
	We assume the operators $\LL$, $\LLd$ admit a complete biorthogonal eigenfunction expansion with respect to the standard inner product $\langle u\given v\rangle=\int_{\R^n}u^*(\mbx)v(\mbx)\,d\mbx$
	\begin{equation}
		\label{eq:spectral-decomposition}
		\LL[P_\lambda]=\lambda P_\lambda,\quad\LLd[Q^*_\lambda]=\lambda Q^*_\lambda,\quad  \langle Q_\lambda\given P_{\lambda'}\rangle=\delta_{\lambda\lambda'}.
	\end{equation}
	The eigenmode $P_0$ represents the unique stationary probability distribution. 
		The normalization eq.~\eqref{eq:spectral-decomposition} implies $Q_0\equiv 1$.  
	The eigenfunction expansion allows us to write the density exactly as 
	\begin{equation}\label{eq:condDensity}
		\rho(\textbf{y},t|\textbf{x},s) = P_0(\textbf{y}) + \sum_{\lambda\not=0} e^{\lambda(t-s)} P_\lambda(\textbf{y}) Q^*_\lambda(\textbf{x}).
	\end{equation}
	
	As in \cite{thomas2014asymptotic} we assume that the system is ``robustly oscillatory", meaning (i) the nontrivial eigenvalue with least negative real part $\lambda_\pm = \mu \pm i\omega$ is complex (with $\omega > 0$), (ii) $|\omega/\mu| \gg 1$, and (iii) for all other eigenvalues $\lambda'$, $\Re[\lambda'] \leq 2\mu$. 
	These conditions guarantee that we can extract the ``stochastic asymptotic phase" $\psi(\mbx)$ from the eigenfunctions $Q^*_{\pm}$ corresponding to
	$\lambda_\pm$ by writing $Q^*_\pm(\mbx)=|Q^*_\pm|e^{\pm i\psi(\mbx)}$, so $\psi(\mbx) = \pm \arg(Q^*_\pm(\mbx))$. 
	Along trajectories, the eigenfunctions $Q^*_\pm(\mbX(t))$  evolve in the mean  as
	$\frac{d}{dt}\E[Q_\pm^*] = \lambda_\pm\E[Q_\pm^*]$, so they behave as a linear focus, oscillating with a period of $2\pi/\omega$ and decaying as $1/|\mu|$ as the density approaches the steady state $P_0$.
	
	We now turn to the special case of planar systems \footnote{In the $n>2$ case we expect $n-1$ stochastic amplitudes $\Sigma^i$ corresponding to the $n-1$ deterministic Floquet modes. 
		The effective vector field would be determined by a system of $n$ equations, e.g.~$\nabla Q^*_\pm(\textbf{x}) \cdot \mbF(\textbf{x}) = \lambda_\pm Q^*_\pm(\textbf{x}), \enskip   \nabla \Sigma^i(\textbf{x}) \cdot \mbF(\textbf{x}) = \lambda^i_\text{Floq} \Sigma^i(\textbf{x})$, for $1\le i \le n-1$.  
		If $|\lambda_\text{Floq}^1|\ll|\lambda_\text{Floq}^i|$ for $i>1$ we expect the flow generated by $\mbF$ will have a 2D invariant manifold $\Sigma^i=0, i>1$, on which the dynamics will be well approximated by our 2D (phase, amplitude) construction.}. 
	In order to establish the amplitude part of the phase-amplitude reduction, we require the slowest mode describing pure contraction without an associated oscillation. 
		In analogy with \cite{thomas2014asymptotic} we seek the \textit{asymptotic} amplitude behaviour.
	Assuming that $\LLd$ has a unique, real, least negative eigenvalue, $\lambda_\text{Floq}\leq2\mu<0$, 
	its corresponding eigenfunction, denoted as  $\Sigma(\mbx)\equiv Q_{\text{Floq}}(\mbx)$, will decay in the mean as 
	\begin{equation}
		\frac{d}{dt}\E[ \Sigma]=\lambda_\text{Floq}\E[\Sigma],
	\end{equation}
	cf. eq.~\eqref{eq:phase-amplitude-deterministic}.
	We interpret $\Sigma(\mbx)$ as the generalization of the amplitude (or isostable) coordinate for the stochastic system eq.~\eqref{eq:SDE}.
	If $\lambda_2<\lambda_\text{Floq}$ is a more negative real eigenvalue, then the isostable mode dominates the non-oscillatory convergence on time scales $\tau\gtrapprox 1/(\lambda_\text{Floq}-\lambda_2)$.
	
	Moreover, the average of individual trajectories, when transformed to the amplitude variable, shows a purely exponential decay towards the set $\Sigma_0\equiv \{\mbx\given\Sigma(\mbx)=0\}$.
	%Thus the set $\Sigma_0$ 	describes the  mean amplitude behavior of an ensemble of trajectories, over the long term.
	Since $\Sigma_0$ in the deterministic case corresponds to a LC, it is natural to ask what $\Sigma_0$ represents in the stochastic system. To that end, we define a vector field $\mathbf{F}$ via
	\small
	\begin{equation}
		\label{eq:effective-vector-field-condition}
		\nabla Q^*_\pm(\textbf{x}) \cdot \mbF(\textbf{x}) = \lambda_\pm Q^*_\pm(\textbf{x}), \enskip   \nabla \Sigma(\textbf{x}) \cdot \mbF(\textbf{x}) = \lambda_\text{Floq} \Sigma(\textbf{x}).
	\end{equation} \normalsize
	Eq.~\eqref{eq:effective-vector-field-condition} determines a unique vector field, provided the gradients of $\Sigma(\mbx)$ and $Q^*_\pm(\mbx)$ are linearly independent at each point $\mbx$. 
		Although the original deterministic system $\dot{\mbx}=\mbf(\mbx)$ may not have oscillatory dynamics, if the full system eq.~\eqref{eq:SDE} is robustly oscillatory, the deterministic vector field $\Re[\mbF]$ will generate a flow with a stable LC, coinciding with $\Sigma_0$.  
		Moreover, eq.~\eqref{eq:effective-vector-field-condition} generates the same amplitude dynamics as $\E[\Sigma]$: pure exponential decay, at rate $\lambda_\text{Floq}$, towards $\Sigma_0$. %We denote $\Re[\mbF]$ as the effective vector field since it shows how the noise effectively modifies the underlying dynamics of the system to sustain oscillations around the effective limit cycle $\Sigma_0$.	% cut: (argument here that the columns of the matrix (nabla theta,nabla sigma) are linearly independent, so we can invert (matrix x F=vector) to solve for F.) \textcolor{blue}{Above equation is correct for the MRT phase. TL phase will satisfy $\mbF(\textbf{x})\cdot  \nabla \psi(\textbf{x}) = \omega + D\langle\nabla \log(u(\mbx))|\nabla \psi(\mbx)\rangle$ with $u(\mbx) = |Q^*_\lambda(\mbx)|$.}
		In the remainder of the paper (see SI for numerical details) we illustrate the construction of $\Sigma$, $Q^*_\pm$, $\Re[\mbF]$ and $\Sigma_0$ for a system arising from an underlying LC as well as for noise-induced oscillations and show that it has properties comporting well with physical intuition. Indeed, in the special case of a LC system perturbed by noise, we observe that $\Sigma$ converges to $\sigma$ as $D\to 0$, in each example we have studied.
	
	\begin{figure}[t!]
		\begin{center}
			\includegraphics[width=9cm]{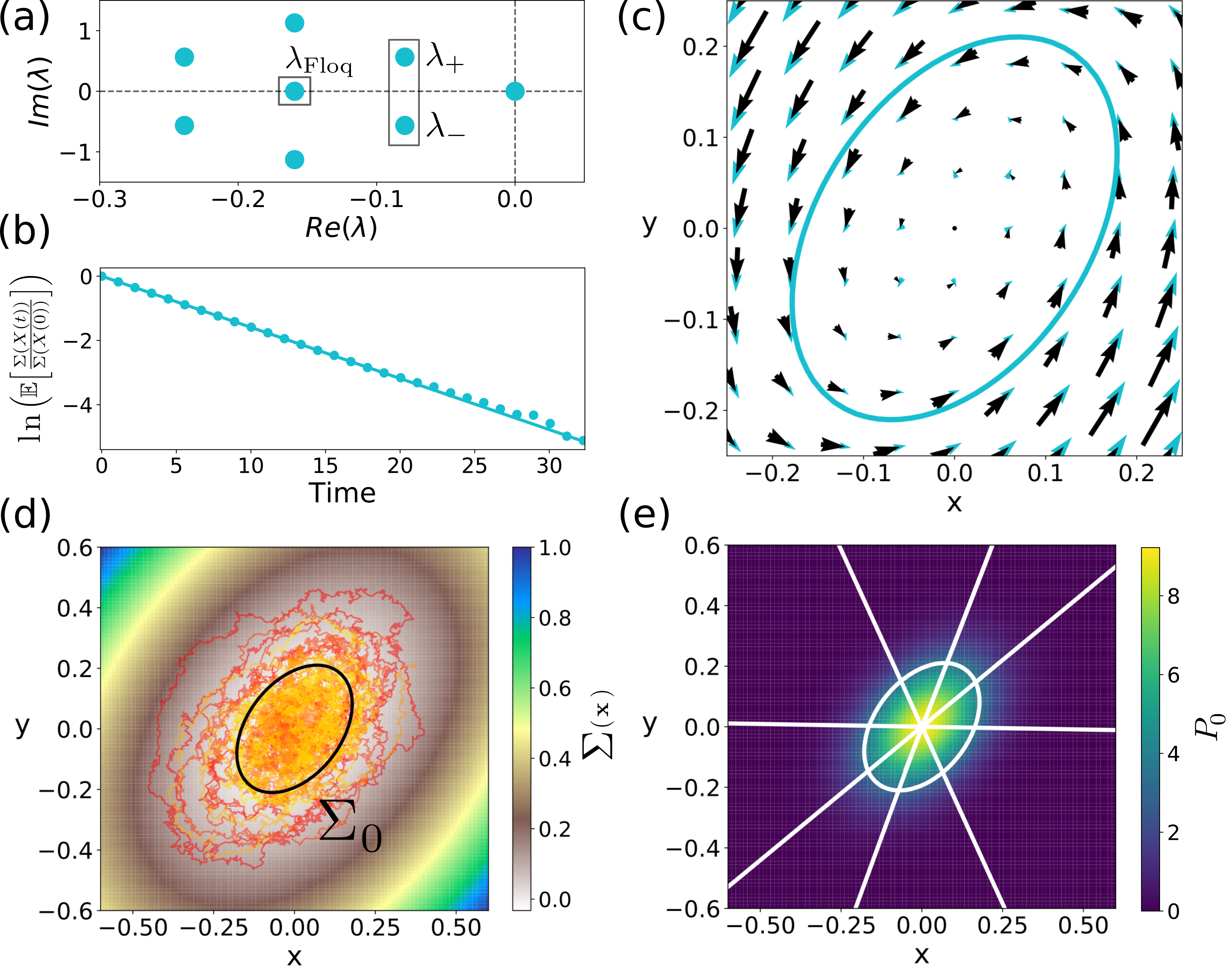}
			\caption{Noisy linear focus eq.~\eqref{eq:langEq} with coefficients  $A=[0.1598, -0.52; 0.7227, -0.319]$ (taken from \citep{powanwe2019determinants}), and isotropic noise $B=\sqrt{2D}\cdot[1,0;0,1]$. 
				$D=1.25\times 10^{-3}$ for all panels.
				(a) Eigenvalue spectrum of  $\mathcal{L}^\dagger$. 
				$\lambda_\pm$ and $\lambda_{\text{Floq}}$ are the complex and real eigenvalues having least negative real part.
				%The eigenfunctions associated to the slowest complex and real eigenvalues allow to define a phase $\Theta(\textbf{x})$ and isostable $\Sigma(\textbf{x})$ function. 
				(b) $\ln\left(\E\left[\Sigma(\mbX(t))/\Sigma(\mbX(0))\right]\right)$ versus time, showing exponential decay of the mean isostable coordinate.
				(c) Deterministic vector field (black) and effective vector field $\mbF(\mbx)$ and its associated effective limit cycle (blue).
				(d) Isostable function and ten  trajectories. 
				Black oval marks the level curve $\Sigma_0$. 
				(e) Stationary probability distribution (color coded), eight isochrons (straight lines) and $\Sigma_0$ (closed line).}
			\label{fig:focus}
		\end{center}
	\end{figure}
	
	\textit{Phase--Amplitude description of a spiral sink -} As a first illustration of the isostable construction for a stochastic oscillator we study a two-dimensional Ornstein-Uhlenbeck process with complex eigenvalues at the origin. This system would not oscillate in the absence of noise.
	Surprisingly, however, phase reduction via the spectral decomposition eq.~\eqref{eq:spectral-decomposition} is well defined \cite{thomas2019phase}.  
	
	We consider a Langevin equation in the form:
	\begin{equation}\label{eq:langEq}
		\dot{\mbX} = A\mbX + B\xi(t)
	\end{equation}
	and assume that the two eigenvalues of $A$ form a complex conjugate pair, $\lambda_\pm = \mu \pm i \omega$ (cf.~Fig.~\ref{fig:focus}a)
	\begin{equation}
		A = \begin{pmatrix} \mu & -\omega \\ \omega & \mu \end{pmatrix}
		, \quad B = \begin{pmatrix} B_{11} & B_{12} \\ B_{21} & B_{22} \end{pmatrix}
	\end{equation} 
	with $B_{ij}$ arbitrary. The stationary probability density $P_0$ of eq.~\eqref{eq:langEq} has a planar Gaussian shape (Fig.~\ref{fig:focus}e) \cite{gardiner1985handbook}. As the spectrum of $\LLd$ shows (Fig.~\ref{fig:focus}a) the  first complex mode is part of an entire family of harmonics \cite{leen2016eigenfunctions}. By using the eigenfunctions $Q^*_\pm(\textbf{x}) = \mbx_1 \pm i \mbx_2$ associated to the smallest complex eigenvalue $\lambda_\pm = \mu \pm i \omega$, we can define an asymptotic phase function $\psi(\textbf{x}) = \arctan(\mbx_2/\mbx_1)$  (Fig.~\ref{fig:focus}e) \cite{thomas2019phase}. In addition, we find the eigenfunction $\Sigma(\textbf{x})$ of $\mathcal{L}^\dagger$, associated with the least negative nontrivial \textit{real} eigenvalue $\lambda_{\text{Floq}}$, to be
	\begin{equation}
		\Sigma(\mbx) = 2 + \frac{\mu}{\epsilon}(\mbx^2_1 + \mbx^2_2), \qquad
		\lambda_{\text{Floq}} = 2\mu,
	\end{equation}
	depicted in Fig.~\ref{fig:focus}d, with $\epsilon = (B^2_{11} + B^2_{12} + B^2_{21} + B^2_{22})/4$. Hence, using $Q_\pm^*$ and $\Sigma$
	we can obtain an analytical expression for the effective vector field $\mbF(\mbx)$ in eq.~\eqref{eq:effective-vector-field-condition}, \small
	\begin{equation}
		\begin{aligned}\label{eq:eff_vect_fieeld_OU}
			\mbF_{1}(\textbf{x}) &= \mu \mbx_1 - \omega \mbx_2 + 2\epsilon\frac{\mbx_1-i \mbx_2}{\mbx^2_1 + \mbx^2_2}, \\
			\mbF_{2}(\textbf{x}) &= \omega \mbx_1 + \mu \mbx_2 +  2\epsilon\frac{\mbx_2+i \mbx_1}{\mbx^2_1 + \mbx^2_2}.
		\end{aligned}
	\end{equation} \normalsize
	The real part of $\mbF$, plotted in Fig.~\ref{fig:focus}c, shows how the dissipative effect of the dynamics ($\mu<0$) combines with the expansive effects of the noise ($\epsilon>0$) to give a finite effective radius, 
	(which is evident upon writing eq.~\eqref{eq:eff_vect_fieeld_OU} in polar coordinates) 
	coinciding with the zero level curve of the isostable function $\Sigma(\mbx)$. 
	That is, 
	\begin{equation}
		\Sigma_0(\mbx) =  \mbx^2_1 + \mbx^2_2 = \frac{2\epsilon}{|\mu|}
	\end{equation} 
	with Floquet exponent $\lambda = 2\mu = \lambda_{\text{Floq}}$. 
		Hence, averaging an amplitude using $\Sigma(\mbX(t))$ instead of the original variables $\mbX(t)$, provides a meaningful asymptotic amplitude $\Sigma_0$ (the effective limit cycle) to which trajectories decay exponentially in the mean (see Fig.~\ref{fig:focus}b,d).

	%whereas averaging an amplitude in terms of the original variables is not meaningful averaging trajectories in terms of the X gives a meaningful amplitude average whereas the average of trajectories in terms of the original variables would in the long run lead to the origin, 

	%Hence, we see how from the This set corresponds to the stochastic limit cycle for eq.~\eqref{eq:langEq} (see Fig.~\ref{fig:focus}e). Stochastic trajectories converge exponentially in the mean to $\Sigma_0$ (Fig.~\ref{fig:focus}b).
	
	%Fig.~\ref{fig:focus} illustrates the elements of our phase-amplitude reduction for the spiral focus \eqref{eq:langEq}.

	%\textcolor{blue}{Discuss Fig 1} Spiral-sink systems with noise-sustained oscillations, also known as quasicycles appear at multiple fields .

\begin{figure}[t]
		\begin{center}
			\includegraphics[width=9cm]{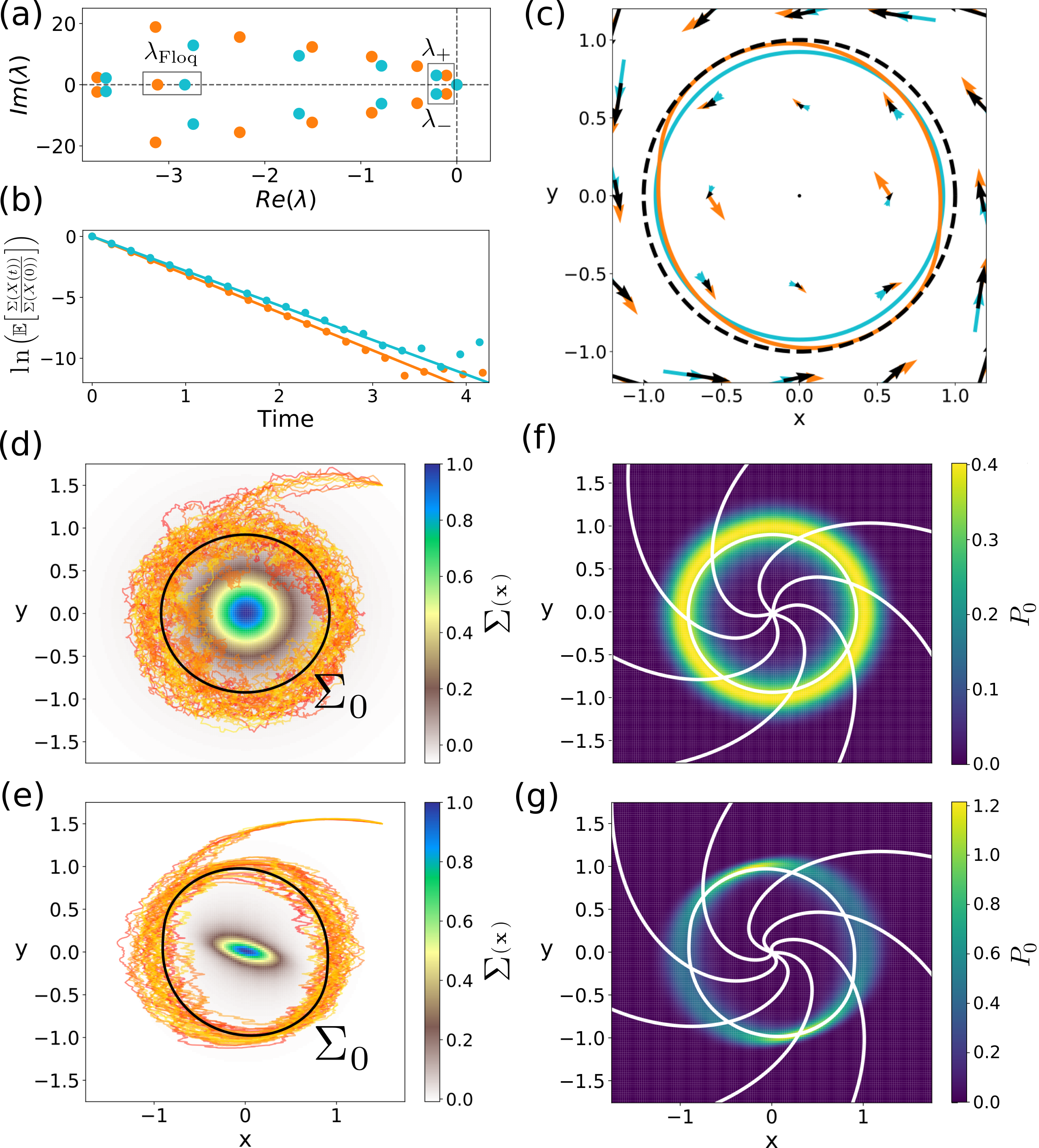}
			\caption{Stuart-Landau oscillator with isotropic noise ($D_{x,y}=0.1$, blue symbols) and anisotropic noise ($D_x=0.1$, $D_y=2.5\times 10^{-4}$, orange symbols). (a) Eigenvalue spectra of $\mathcal{L}^\dagger$ for isotropic and anisotropic noise. (b) $\ln\left(\E\left[\Sigma(\mbX(t))/\Sigma(\mbX(0))\right]\right)$ versus time, showing exponential decay of the mean isostable coordinate. (c) Deterministic limit cycle and vector field (black) and effective vector fields $\mbF$ and their associated effective limit cycles (d-e) Isostable functions and ten trajectories for isotropic (d) and anisotropic (e) noise. (f-g) Stationary probability distribution (color coded), with seven isochrons (straight lines) and $\Sigma_0$ (closed line) for isotropic (f) and anisotropic (g) noise.}
			\label{fig:slPanel}
		\end{center}
	\end{figure}

	\textit{Noisy Stuart-Landau (SL) Oscillator -}
	The SL equations capture universal dynamics near a Hopf bifurcation,
%	\rot{which arise as an example of coherence resonance \cite{ushakov2005coherence}}
	%
	\begin{equation*}
		%\label{eq:SL}
		\small
		\begin{aligned}
			\dot{\mbX} &= b \mbX (1 - (\mbX^2+\mbY^2)) - \mbY(1 + b a (\mbX^2+\mbY^2))+\sqrt{2D_x}\xi_x(t),\\
			\dot{\mbY} &= b \mbY (1 - (\mbX^2+\mbY^2)) + \mbX(1 + b a (\mbX^2+\mbY^2)) +\sqrt{2D_y}\xi_y(t)
		\end{aligned}
		\normalsize
	\end{equation*}
	with $a, b = [1, 2] \in \mathbb{R}$. In the absence of noise this system %\eqref{eq:SL} 
	has a LC of period $T = 2 \pi/(1 + ba)$ and a Floquet exponent $\lambda = -2b$. 
	To see how isostables capture the effects of noise, in Fig.~\ref{fig:slPanel} we study an isotropic ($D_{x,y}=0.1$) and an anisotropic ($D_x=0.1$, $D_y=2.5\times 10^{-4}$) noise case.
	
	%\textcolor{blue}{Discuss Fig.2}
	%\apc{Acaba de completar bien los casos del SL y del focus (con las referencias) y pensar el calculo teorico del caso del focus} \apc{Piensa como modificar la discussion para añadir lo que falta}
	
	The isotropic case has larger total noise than the anisotropic case. For this reason, for the anisotropic case, $\Sigma_0$ and $\lambda_\text{Floq}$ lie nearer their deterministic analogues than for the isotropic case (Fig.~\ref{fig:slPanel}(a,c)). Under reduced total noise, the real part of the smallest complex eigenvalue pair increases, indicating slower decay of coherent oscillation to the steady-state distribution. We notice the occurrence of a family of eigenvalues $\lambda_k \approx \pm i\omega k - \mu k^2$  \cite{thomas2014asymptotic, vcrnjaric2019koopman}.
	At the same time, $\lambda_\text{Floq}$ becomes \emph{more negative}, indicating a faster decay of the isostable coordinate (cf.~Fig.~\ref{fig:slPanel}(a,b)).  
	As in the spiral focus system, the slowest decaying real and complex eigenfunctions of the noisy SL system lead to an effective vector field $\mbF$ that exhibits a stable limit cycle, $\Sigma_0$ (panel c). 
	Note the rotational symmetry of $\Sigma_0$ in the isotropic case (blue circle in c) versus the lack of symmetry of the anisotropic case (orange ellipse in c); the dashed black curve shows the deterministic SL limit cycle, for comparison.  
	In both cases, the noise reduces the effective radius of the LC.
	The asymmetry of the anisotropic case appears in all level sets of $\Sigma$, as well as the stationary distribution and isochrons (panel e,g).

\begin{figure}[t]
		\begin{center}
			\includegraphics[width=9cm]{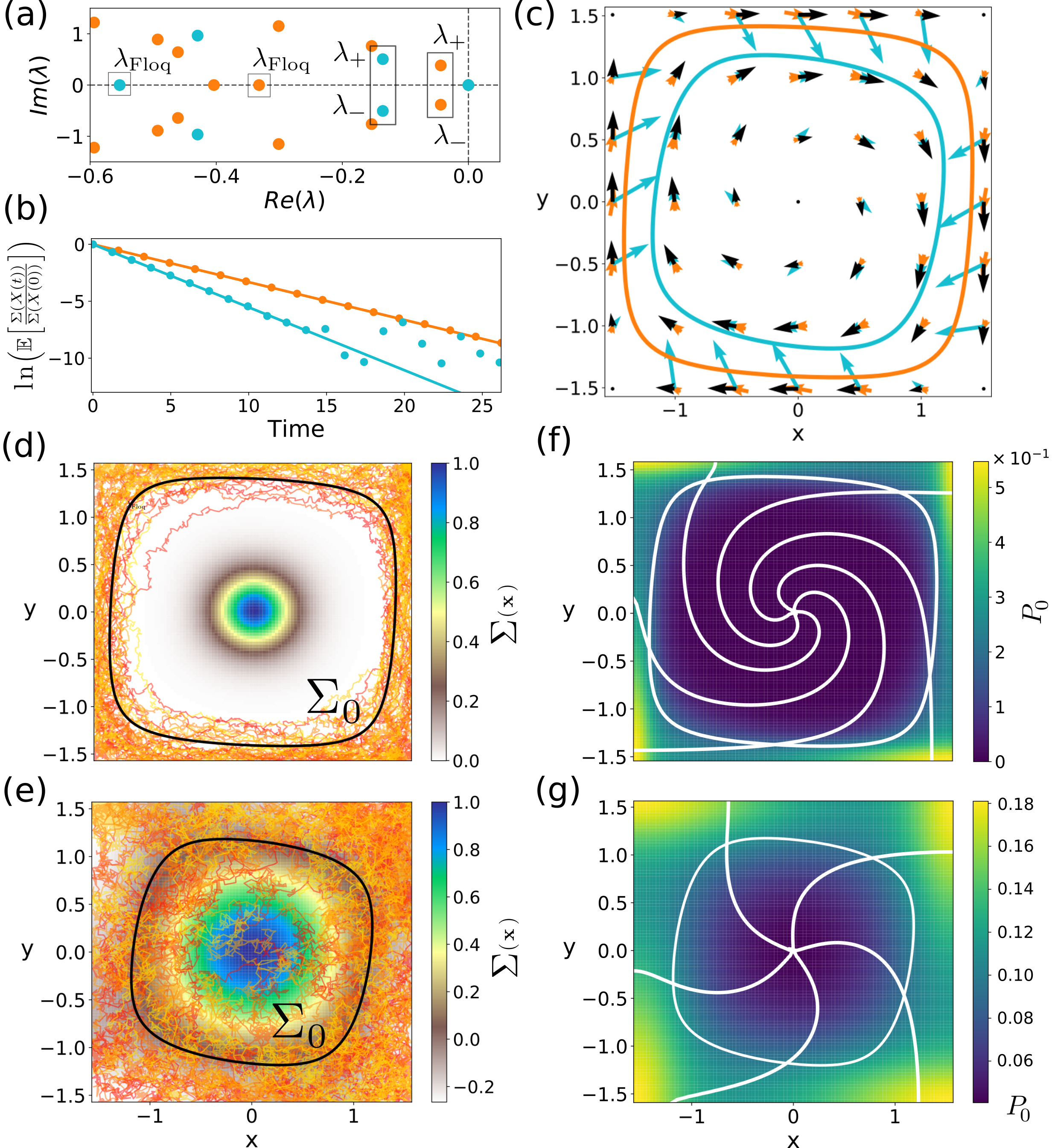}
			\caption{Heteroclinic oscillator with lower ($D=0.01125$, orange symbols) and higher ($D=0.1$, blue symbols) noise. (a) Eigenvalue spectra of $\mathcal{L}^\dagger$ for lower and higher noise levels. 
				In this case, reducing $D$ shifts both $\Re[\lambda_\pm]$ and $\lambda_\text{Floq}$ to the right. %The eigenfunctions associated to the slowest complex and real eigenvalues allow to define a phase $\Theta(\textbf{x})$ and isostable $\Sigma(\textbf{x})$ function.  
				(b) $\ln\left(\E\left[\Sigma(\mbX(t))/\Sigma(\mbX(0))\right]\right)$ versus time, showing exponential decay of the mean isostable coordinate. (c) Deterministic vector field $\mbf(\mbx)$ (black) and effective vector fields $\mbF$ and their associated effective limit cycles. $\mbF$ and $\mbf$ differ, especially near the borders. (d-e) Isostable functions and ten trajectories for lower (d) and higher (e) noise. 
				As noise increases, trajectories spread across the phase space, and $\Sigma_0$ moves away from the boundaries. 
				(f-g) Stationary probability distribution (color coded), five isochrons (straight lines) and $\Sigma_0$ (closed line) for lower (f) and higher (g) noise level. %Note the increased isochron curvature at lower noise.
			}
			\label{fig:hetPanel}
		\end{center}
	\end{figure}

	\textit{Noisy Heteroclinic Oscillator -} The underlying vector field of a heteroclinic oscillator creates a closed loop of trajectories connecting a sequence of saddle equilibria \cite{Afraimovich2004heteroclinic,
		Armbruster2003noisy,
		MayLeonard1975nonlinear}.
	Fig.~\ref{fig:hetPanel} shows the phase-amplitude analysis of the heteroclinic oscillator
	system 
	\begin{equation} 
		\label{eq:noisy-het}
		\begin{aligned}
			\dot{\mbX} &=  \cos(\mbX)\sin(\mbY) + \alpha \sin(2 \mbX) + \sqrt{2 D}\xi_1(t) \\
			\dot{\mbY} &= -\sin(\mbX)\cos(\mbY) + \alpha \sin(2 \mbY) + \sqrt{2 D}\xi_2(t) 
		\end{aligned}
	\end{equation} 
	with $\alpha = 0.1$ and reflecting boundary conditions on the domain $-\pi/2 \leq \{\mbX, \mbY\} \leq \pi/2$. 
	Without noise the system has an attracting heteroclinic cycle. 
	Therefore, it can only sustain robust oscillations in the presence of noise \cite{thomas2014asymptotic}. 
	We study this oscillator for a lower ($D = 0.01125$) and for a higher level of noise ($D = 0.1$). 
	For the lower level of noise we observe a rightward shift of both the smallest complex eigenvalue, indicating a longer persistence time of the oscillation, and also the smallest real eigenvalue, indicating slower convergence to the stationary state and slower decay of the isostable coordinate towards $\Sigma_0$ (cf.~\ref{fig:hetPanel}(a,b)) \footnote{In \cite{thomas2014asymptotic} Thomas and Lindner numerically solved the eigenvalue problem eq.~\eqref{eq:spectral-decomposition} for the heteroclinic system eq.~\eqref{eq:noisy-het} using a Fourier mode decomposition method.  
		%Care must be taken to implement the correct boundary conditions.  
		Due to a subtle error in our treatment of the boundary conditions, the slowest decaying \emph{real} eigenvalue plotted in Fig.~2(c) of \cite{thomas2014asymptotic} was incorrect. It has been corrected in Fig.~\ref{fig:hetPanel}(a) of this Letter. This error had no effect on the analysis or conclusions in \cite{thomas2014asymptotic}, which concerned only the complex-valued eigenvalue and its eigenfunction.}. At the same time, the mean period increases logarithmically in $D$ (not shown) \cite{Giner-Baldo2017JStatPhys}.
	The effective vector fields $\mbF(\mbx)$  differ significantly from the deterministic system $\mbf(\mbx)$ (panel c), especially near the domain boundaries. 
	Noise causes the effective direction of flow to change, from  paralleling the walls, to pointing inwards towards the origin.
	Consequently, although $\mbf$ does not support a limit cycle, $\mbF$ does.
	
	For smaller noise, the stationary distribution (panel f) is pressed tightly against the domain boundaries, and trajectories rarely visit the interior of the domain (panel d). Hence, the zeroth isostable $\Sigma_0$ remains close to the domain walls (panel d).  
	For larger noise, the stationary distribution spreads farther from the walls. Due to the reflecting boundaries, trajectories visit all regions of the domain, and $\Sigma_0$ contracts partway into the domain interior (panel e).
	%As the level curves of $\Sigma$ show, as the noise increases, trajectories evolve across the whole phase space. Whereas for lower noise trajectories are distributed near the borders, increasing the noise approaches trajectories to the origin. As a consequence, the $\Sigma_0$ curve shrinks as noise amplitude increases (panel d,e). 
	Increasing the noise amplitude reduces significantly the curvature of the isochrons (panel f, g).
	
	\textit{Discussion -} 
	In this Letter we advance the work begun in \cite{thomas2014asymptotic} towards a self-consistent phase-amplitude description of stochastic oscillations.   Our methodology applies to stochastic oscillations, even if they are noise induced. 
	By providing an analog of both phase and amplitude, we establish a unified description of deterministic and stochastic oscillations. 
	The zeroth isostable $\Sigma_0$ links the deterministic setting, in which the amplitude coordinate obeys $\dot{\sigma}\propto -\sigma$, and the stochastic setting, in which $\dot{\langle\Sigma\rangle}\propto - \langle\Sigma\rangle$. 
	Thus, for an oscillatory stochastic system, trajectories asymptotically approach $\Sigma_0$  ``in the mean".
	%Of particular relevance in this letter is the zero-th isostable $\Sigma_0$. 
	In addition to unifying noise-perturbed deterministic LCs and noise-induced oscillations in heteroclinic systems, $\Sigma_0$ also captures the mean amplitude dynamics in  
	the case of quasicycles \cite{LugoMcKane2008PREquasicycles,BrooksBressloff2015quasicycles,Bressloff2010PREquasicycles,duchet2021optimizing}. 
	Together with \cite{thomas2014asymptotic}, the framework introduced in this Letter may provide a basis for studying noise-driven oscillations in excitable systems \cite{lindner2004effects}, stochastic bifurcations \cite{arnold1995random}, coherence resonance \cite{PikKur97,ushakov2005coherence}, or to construct a physically grounded definition for stochastic limit cycles.

	Our results assume a particular structure for the spectrum of $\LLd$: the ``robustly oscillatory criteria'' are met and the initial slowest decaying complex pair entails an entire family of eigenvalues \cite{thomas2014asymptotic}. After this family constituting the phase dynamics, then, the slowest decaying part of the remaining eigenmode decomposition is precisely the least negative real eigenvalue. 
		Moreover, the structure of the spectrum of $\LLd$ can  also be interpreted in terms of the deterministic phase-amplitude coordinates derived from the eigenfunctions of the Koopman operator \cite{mauroy2018global, shirasaka2020phase}. 
	Indeed, the connection between the asymptotic phase for stochastic oscillators introduced in \cite{thomas2014asymptotic} and the Koopman operator has been noted  \cite{KatoNakao2020quantum,EngelKuehn2021CommMatPhys,kato2021asymptotic}.
	This connection is not coincidental: the backward Kolmogorov operator at the heart of our analysis, $\mathcal{L}^\dagger$, is also the generator of the Markov process, as well as the stochastic Koopman operator \cite{vcrnjaric2019koopman}. 
	Therefore, beyond its theoretical interest, since our phase-amplitude functions are encoded in the Koopman operator, modern methods for extracting Koopman eigenfunctions from data may lead to practical methods for establishing phase-amplitude coordinates for noisy oscillatory systems in medicine, biology, engineering, economics, control, or other areas \cite{budivsic2012applied, mauroy2020koopman,proctor2016dynamic,schmid2010dynamic}.
	%our approach may have interesting implications in many real world applications such as data-science, medical or climate among others fields in which stochastic oscillations are present. 
	%Indeed, during the last decade there is a growing interest in Koopman operator, being nowadays an indispensable tool to analyse data and control dynamics \cite{budivsic2012applied, mauroy2020koopman}. While our approach can benefit from the current numerical methodology, the dynamical mode decomposition \cite{proctor2016dynamic, schmid2010dynamic} in data, thus providing a theoretical framework allowing for a complementary interpretation of results, is also an appealing topic for further applied research.
	
	%\begin{acknowledgments}
		%\textit{Acknowledgments -} APC acknowledges support from the Basic Research Program of the National Research University Higher School of Economics as well as IFISC, UIB-CSIC for hosting him.This work was supported in part by NSF grant DMS-2052109 and by NIH BRAIN Initiative grant R01-NS118606. PT thanks the Oberlin College Department of Mathematics for research support.
		
		\textit{Acknowledgments -} Basic Research Program of the National Research University Higher School of Economics and IFISC, UIB-CSIC (APC). NSF grant DMS-2052109, NIH BRAIN Initiative grant R01-NS118606 and the Oberlin College Department of Mathematics (PT).

		%\end{acknowledgments}
	
	% Produces the bibliography via BibTeX.

%\bibliography{apssamp}

%merlin.mbs apsrev4-1.bst 2010-07-25 4.21a (PWD, AO, DPC) hacked
%Control: key (0)
%Control: author (8) initials jnrlst
%Control: editor formatted (1) identically to author
%Control: production of article title (-1) disabled
%Control: page (0) single
%Control: year (1) truncated
%Control: production of eprint (0) enabled
%

\appendix

\widetext

\section*{Supplementary Material \\ Numerical Results}

To obtain the numerical results in the main manuscript, we followed the procedure introduced in  \cite{cao2020partial}. Given a Langevin equation as in eq.~(3) in the main manuscript, we first chose a rectangular domain
\begin{equation}\label{eq:domainD}
	\mathcal{D} = [ x^-, x^+ ] \times [y^-, y^+].   
\end{equation} 
For the spiral sink and the Stuart-Landau oscillator, in which the phase space is unbounded, we use a truncated domain whose size is chosen large enough so that the probability for individual trajectories $\mbX(t)$ to reach the boundaries is very low. For the special case of the heteroclinic oscillator boundaries are given by the nature of the system. 

By discretizing the domain $\mathcal{D}$ in $N$ and $M$ points such that $\Delta x = (x_+ - x_-)/N$ and $\Delta y = (y_+ - y_-)/M$, we can build $\LLd$ (and/or $\LL$) by using a standard finite difference scheme. In general we used centered finite differences, that is
\begin{align}
	(\partial_x T)_{i,j} &= \frac{T_{i+1,j} - T_{i-1,j}}{2\Delta x}, \\ (\partial_{xx} T )_{i,j} &= \frac{T_{i+1,j} - 2T_{i,j} + T_{i-1,j}}{(\Delta x)^2}
\end{align}
except at the borders of the domain. For the heteroclinic oscillator, we implemented adjoint reflecting boundary conditions in the borders of the domain. By contrast, for the unbounded systems, since there is no natural border, we just substituted the centered finite difference scheme by a forward (or backward) finite difference scheme. Using adjoint reflecting boundary conditions for these systems yielded numerically very similar results.

After diagonalizing the resulting $(N \cdot M, N \cdot M)$ matrix, we obtain the eigenvalues and the associated eigenfunctions of $\LLd$ ($\LL$) (for parameter values and resulting eigenvalues see Table~\ref{tab:table1}). We recall that we are not interested in the complete spectrum of $\LLd$ ($\LL$). For $\LLd$ we just consider (and hence present in Figs.~1-3(a)) the part of the spectrum which is relevant for our analysis. That is, we consider the eigenvalue associated with the slowest decaying complex eigenfunction $Q^*_\pm(\mbx)$ (from which we obtain the asymptotic phase $\psi(\mbx) = \pm \arg(Q^*_\pm(\mbx)$) and the eigenfunction $\Sigma(\mbx)$ associated with the least negative purely real eigenvalue $\lambda_{\text{Floq}}$. For $\LL$ we just consider the eigenmode associated with the eigenvalue $\lambda=0$ which gives the stationary probability distribution $P_0$.\\
\begin{table}[h]
\begin{ruledtabular}
	\begin{tabular}{cccccccccc}
		&$N$ &$M$ &$x_+$&
		$x_-$ &$y_+$ &$y_-$ &$\mu$ &$\omega$ &$\lambda_\text{Floq}$\\
		\hline
		Sp. Sink & 151 & 151 & 0.6 & -0.6 & 0.6 & -0.6 & -0.080 & 0.564 & -0.159 \\
		SL-iso & 151 & 151 & 1.75 & -1.75 & 1.75 & -1.75 & -0.213 & 3.032 & -2.833\\
		SL-ani & 151 & 151 & 1.75 & -1.75 & 1.75 & -1.75 & -0.108 & 3.008 & -3.117\\
		Het-low & 151 & 151 & $\pi/2$ & $-\pi/2$ & $\pi/2$ & $-\pi/2$ & -0.044 & 0.383 & -0.332\\
		Het-high & 151 & 151 & $\pi/2$ & $-\pi/2$ & $\pi/2$ & $-\pi/2$ & -0.136 & 0.505 & -0.553\\
	\end{tabular}
	\caption{\label{tab:table1}Numerical values of the parameters and resulting leading eigenvalues for the different stochastic oscillators.}
\end{ruledtabular}
\end{table}

%\end{widetext}

\end{document}